\documentclass[12pt]{amsart}
\usepackage{amsmath,amscd,amssymb,amsfonts,graphics}
\newcommand{\h}{\hbox}
\newcommand{\q}{\quad}

\newcommand{\bs}{\par\bigskip}
\newcommand{\ms}{\par\medskip}
\newcommand{\sk}{\par\smallskip}
\newcommand{\bsn}{\par\bigskip\noindent}
\newcommand{\msn}{\par\medskip\noindent}

\newcommand{\ssb}{\raise.15ex\h{${\scriptscriptstyle\bullet}$}}
\newcommand{\ssc}{\,\raise.15ex\h{${\scriptstyle\circ}$}\,}
\newcommand{\msum}{\hbox{$\sum$}}

\newcommand{\mopl}{\hbox{$\bigoplus$}}
\newcommand{\C}{{\mathbb C}}
\newcommand{\N}{{\mathbb N}}
\newcommand{\PP}{{\mathbb P}}

\newcommand{\Z}{{\mathbb Z}}
\newcommand{\J}{{\mathcal J}}
\newcommand{\OO}{{\mathcal O}}

\newcommand{\dd}{\partial}
\newcommand{\ddd}{{\rm d}}
\newcommand{\df}{{\rm d}f}
\newcommand{\mm}{{\mathfrak m}}
\newcommand{\muc}{\check{\mu}}
\newcommand{\nuc}{\check{\nu}}
\newcommand{\rhoc}{\check{\rho}}
\newcommand{\chic}{\check{\chi}}
\newcommand{\Euc}{\check{\rm E}{\rm u}}
\newcommand{\reg}{{\rm reg}}
\newcommand{\Der}{{\rm Der}}
\newcommand{\Diff}{{\rm Diff}}
\newcommand{\Gr}{{\rm Gr}}
\newcommand{\Om}{\Omega}

\newcommand{\bl}{\bigl}
\newcommand{\br}{\bigr}

\newcommand{\ges}{\geqslant}
\newcommand{\les}{\leqslant}
\newcommand{\1}{\hskip1pt}
\newcommand{\lc}{{\it loc.\,cit.}}

\begin{document}
\title[Degeneration of pole order spectral sequences]
{Degeneration of pole order spectral sequences\\
for hyperplane arrangements of 4 variables}
\author{Morihiko Saito}
\address{RIMS Kyoto University, Kyoto 606-8502 Japan}
\begin{abstract} For essential reduced hyperplane arrangements of 4 variables, we show that the pole order spectral sequence degenerates almost at $E_2$, and completely at $E_3$, generalizing the 3 variable case where the complete $E_2$-degeneration is known. These degenerations are useful to determine the roots of Bernstein-Sato polynomials supported at the origin. For the proof we improve an estimate of the Castelnuovo-Mumford regularity of logarithmic vector fields which was studied by H.~Derksen and J.~Sidman.
\end{abstract}
\maketitle
\centerline{\bf Introduction}
\bsn
Let $Z\subset\PP^{n-1}$ be a reduced hyperplane arrangement with $n\ges 3$. Let $f$ be a defining polynomial. We have the {\it pole order spectral sequence\1} which is associated with a double complex whose differentials are given by $\ddd$ and $\df\wedge$, see for instance \cite{DiSa}, \cite{nwh}. The abutting filtration is the pole order filtration, which determines the roots of Bernstein-Sato polynomials supported at the origin under some condition, see \cite[Theorem 2]{bfun}.
\sk
In the case $n=3$, it is proved that the spectral sequence degenerates at $E_2$, and the $E_1$-term of the spectral sequence can be obtained from the Hilbert series of the Milnor algebra $R/(\dd f)$ using the self-duality of the $E_1$-term, where $(\dd f)\subset R:=\C[x]$ is the Jacobian ideal generated by the partial derivatives of $f$, see \cite{DiSa}, \cite{wh}. We have moreover the injectivity of the differential $\ddd_1$ except for certain places which are irrelevant to the computation of the roots supported at the origin. These can be used to determine the roots of Bernstein-Sato polynomials except for certain special cases, see \lc
\sk
Assume from now on $n=4$. We can still calculate the $E_1$-term of the spectral sequence quite explicitly and swiftly using computer programs like Macaulay2 or Singular, see \cite{nwh}. However, there is no method to prove the $E_2$-degeneration for this case although there is no counterexample.
In the strongly free divisor case, it is shown in \cite[Theorem 4.10]{nwh} that the spectral sequence degenerates almost at $E_2$, and completely at $E_3$. This was motivated by a conjecture in \cite{DiSt}. In this paper we show that this theorem also holds (with $m=2d{-}1$) in a hyperplane arrangement case; more precisely, we prove the following.
\msn
{\bf Theorem~1.} {\it For essential reduced hyperplane arrangements in $\PP^3$, the pole order spectral sequence degenerates completely at $E_3$ and almost at $E_2\1;$ more precisely, we have}
$$\mu_k^{(2)}=0\,\,\,(k>2d{-}2),\q\nu_k^{(2)}=0\,\,\,(k>3d{-}1),\q\rho_k^{(2)}=0\,\,\,(k>4d{-}2).
\leqno(1)$$
\ms
Here $\mu_k^{(2)}$, $\nu_k^{(2)}$, $\rho_k^{(2)}$ denote the dimensions of certain $E_2$-terms of the spectral sequence, see \cite{nwh}.
The assertion for $\nu^{(2)}_k$ is slightly weaker than (4.3.4), \h{\lc} for $m=2d{-}2$, where the range of $k$ is given by $k>3d{-}2$. This comes from the difficulty to control the $\mu'''_{{\rm def},k}$ in the notation of {\lc}
For the proof of Theorem 1 we use the Castelnuovo-Mumford regularity.
\sk
This work is partially supported by Kakenhi 15K04816.
\sk
In Section 1 we improve an estimate of the regularity of logarithmic vector fields shown in \cite{DeSi}. In Section 2 we prove Theorem 1 using the estimate of the Castelnuovo-Mumford regularity.
\bs\bs
\vbox{
\centerline{\bf 1. Regularity of logarithmic vector fields}
\bsn
In this section we review some basics of Castelnuovo-Mumford regularity, and improve an estimate of the regularity of logarithmic vector fields shown in \cite{DeSi}.}
\msn
{\bf 1.1.~Castelnuovo-Mumford Regularity} (see \cite{Ei}). Let $M$ be a finitely generated graded $R$-module with $R:=\C[x_1,\dots,x_n]$. Take a minimal graded free resolution
$$\to F_l\to\cdots\to F_0\to M\to 0,$$
where $F_j=\mopl_k\,R(-c_{j,k})$ with $c_{j,k}\in\Z$. The {\it Castelnuovo-Mumford regularity\1} is defined by
$$\reg\,M:=\max_{j,k}\,(c_{j,k}-j).
\leqno(1.1.1)$$
\msn
{\bf Remarks~1.1} (i). By the graded version of Nakayama's lemma, the minimality of the graded free resolution is equivalent to the vanishing of the differential of the complex
$$\to\C\otimes_RF_l\to\cdots\to\C\otimes_RF_0\to 0,$$
where $\C=R/\mm\,$ with $\,\mm\subset R$ the graded maximal ideal, see {\lc} This implies that
$${\rm Tor}^R_j(\C,M)=\mopl_k\,\C(-c_{j,k})\q(j\in\N),
\leqno(1.1.2)$$
hence
$$\reg\,M=\max_j\,\bl({\rm max\,deg}\bl({\rm Tor}^R_j(\C,M)\br)-j\br),
\leqno(1.1.3)$$
with ${\rm max\,deg}(E):=\max\bl\{k\in\Z\mid E_k\ne 0\br\}$.
\ms
(ii) We have a formula similar to (1.1.3) using local cohomology (see Theorem~4.3, \lc)
$$\reg\,M=\max_j\,\bl({\rm max\,deg}\bl(H^j_{\mm}M\br)+j\br).
\leqno(1.1.4)$$
Note that $H^j_{\mm}R=0$ for $j\ne n$, and $H^n_{\mm}R=\C[\dd_1,\dots,\dd_n](-n)$. (The assertion (1.1.4) is much more difficult than (1.1.3).)
\msn
{\bf 1.2.~Logarithmic vector fields.} Let $X$ be a smooth complex algebraic variety or a complex manifold, and $D\subset X$ be an effective divisor. We denote by $\Der_X$ the sheaf of algebraic or holomorphic vector fields on $X$. The subsheaf $\Der_X(-\log D)$ of logarithmic vector fields along $D$ can be defined locally by the exact sequence
$$0\to\Der_X(-\log D)\to\Der_X\to\OO_X/\OO_X(-D),
\leqno(1.2.1)$$
where the last morphism is given by $\Der_X\ni\theta\mapsto[\theta(f)]\in\OO_X/\OO_X(-D)$ with $f$ a local defining function of $D$. This subsheaf is {\it independent\1} of a {\it non-reduced\1} structure of $D$ by the Leibniz law for derivations. (This is quite different from the case of logarithmic differential forms.) From now on we assume $D$ {\it reduced}.
\sk
The following is needed to improve the estimate of the Castelnuovo-Mumford regularity of logarithmic vector fields by \cite{DeSi}.
\msn
{\bf Lemma~1.2.} {\it Assume $X=X_1{\times}X_2$ and $D={\rm pr}_1^*D_1+{\rm pr}_2^*D_2$ with ${\rm pr}_i:X\to X_i$ natural projections $(i=1,2)$. Then we have the canonical isomorphism}
$$\Der_X(-\log D)={\rm pr}_1^*\1\Der_{X_1}(-\log D_1)\,\oplus\,{\rm pr}_2^*\1\Der_{X_2}(-\log D_2),
\leqno(1.2.2)$$
\msn
{\it Proof.} Since we have the canonical isomorphism
$$\Der_X={\rm pr}_1^*\1\Der_{X_1}\oplus{\rm pr}_2^*\1\Der_{X_2},
\leqno(1.2.3)$$
any $\theta\in\Der_X$ is written locally as $\theta_1+\theta_2$ with $\theta_i\in{\rm pr}_i^*\1\Der_{X_i}$ ($i=1,2$). We denote also by $f_i$ the pull-back of a local defining function $f_i$ of $D_i$. The condition
$$\theta(f_1f_2)=\theta_1(f_1)f_2+f_1\theta_2(f_2)\,\in\,(f_1f_2)$$
is equivalent to that $\theta_i(f_i)\in(f_i)$ ($i=1,2$), since $\OO_{X,x}$ is a unique factorization domain and $f_1$, $f_2$ are mutually prime. So the isomorphism (1.2.2) follows taking the pull-back of the exact sequence (1.2.1) for $(X_i,D_i)$ by ${\rm pr}_i^*$ ($i=1,2$), since the latter is an exact functor of $\OO$-modules. This finishes the proof of Lemma~(1.2).
\msn
{\bf 1.3.~Regularity of logarithmic vector fields.} Let $D$ be a hyperplane arrangement in $X=\C^n$ (as an algebraic subvariety). Assume $D$ central, that is, $f$ homogeneous. Set $R=\C[x_1,\dots,x_n]$ with $x_1,\dots,x_n$ the coordinates of $X$, and
$$\Der_R(-\log D):=\Gamma\bl(X,\Der_X(-\log D)\br).$$
This is the graded $R$-module of algebraic logarithmic vector fields on $X$, where the grading is given by $\deg x_i=1$ and moreover $\deg \dd/\dd x_i=-1$ in this paper. So {\it the degree $($as well as the regularity$)$ is shifted by} 1, compared with \cite{DeSi,Sc}.
\sk
The following is due to \cite{Sc} in the case $n=3$, and is shown essentially in \cite{DeSi}.
\msn
{\bf Proposition~1.3.} {\it For an essential central reduced hyperplane arrangement in $X=\C^n$, we have the inequality}
$$\reg\,\Der_R(-\log D)\les\deg D-n.
\leqno(1.3.1)$$
\msn
{\it Proof.} This can be proved by induction on $n$ and $\deg D$ using \cite[Corollary 3.7]{DeSi} where each hyperplane is deleted for induction. If we get a {\it non-essential\1} one by deleting a hyperplane from an arrangement, then the latter is defined by $f=x_ng$ with $g\in R':=\C[x_1,\dots,x_{n{-}1}]$ replacing the coordinates $x_1,\dots,x_n$ appropriately, where the deleted hyperplane is defined by $x_n=0$. (Recall that an arrangement which cannot be defined by a polynomial of fewer variables is called essential.) In this case Lemma~(1.2) implies the canonical isomorphism
$$\Der_R(-\log D)=\Der_{R'}(-\log D'){\otimes_{R'}}R\,\oplus\,\Der_{\C[x_n]}(-\log 0){\otimes_{\C[x_n]}}R,
\leqno(1.3.2)$$
with $D':=g^{-1}(0)\subset\C^{n-1}$. The regularity of the right-hand side of (1.3.2) is bounded by
$$\deg D'-(n{-}1)=\deg D-n$$
by inductive hypothesis (where $\deg D\ges n$ since $D$ is essential). So we can use (1.3.2) instead of Corollary 3.7 for this case, and proceed by induction on $n$ and $\deg D$. This finishes the proof of Proposition~(1.3).
\msn
{\bf 1.4. Regularity of Milnor algebras.} In the notation of (1.3), we have a well-known decomposition
$$\Der_R(-\log D)=\Der_R(-\log D)^0\oplus R\1\theta_0.
\leqno(1.4.1)$$
Here $\theta_0:=\msum_{i=1}^n\,\tfrac{1}{d}\1x_i\dd/\dd x_i$ so that $\theta_0(f)=f$, and
$$\Der_R(-\log D)^0:=\{\theta\in\Der_R\mid\theta(f)=0\}.$$
\sk
Set
$$A_f^p:={\rm Ker}(\df\wedge:\Om^p\to\Om^{p+1}),$$
with $\Om^p$ the algebraic differential $p$-forms on $\C^n$. Put $d:=\deg D=\deg f$. We have the isomorphism of graded $R$-modules
$$\Der_R(-\log D)^0=A_f^{n-1}(n),
\leqno(1.4.2)$$
since both are identified (up to a shift of grading) with
$$\bl\{\,g_i\in R\,\,(i\in[1,n])\mid\msum_i\,g_i\,\dd f/\dd x_i=0\br\}.$$
By definition there are two short exact sequence of graded $R$-modules
$$\aligned&\q\q 0\to\df{\wedge}\Om^{n-1}\to\Om^n\to M\to 0,\\ &0\to A_f^{n-1}\to\Om^{n-1}\buildrel{\df\wedge\,\,}\over\longrightarrow\df{\wedge}\Om^{n-1}(d)\to 0,\endaligned
\leqno(1.4.3)$$
where $M$ is as in \cite{nwh}, and is isomorphic to the shifted Milnor algebra $(R/(\dd f))(-n)$, see also \cite{DIM}.
Note that the $\Om^p$ are free $R$-modules with $\reg\,\Om^p=p$ ($p=n{-}1,n$). There are no $\C$-linear relations among the $\dd f/\dd x_i$ by the essentiality of $D$. So the above two short exact sequences can be extended to a minimal graded free resolution of $M$, and we get
$$\reg\,M=\reg\,\df{\wedge}\Om^{n-1}-1=\reg\,A^{n-1}_f(-d)-2,
\leqno(1.4.4)$$
since we have by (1.4.2)
$$\reg\,A^{n-1}_f(-d)=\reg\,\Der_R(-\log D)^0(-d-n)\ges d+n.$$
\sk
Proposition~(1.3) then implies the following.
\msn
{\bf Corollary~1.4.} {\it For an essential central reduced hyperplane arrangement $D$ in $\C^n$, we have}
$$\reg\,M\les\,2d-2.
\leqno(1.4.5)$$
\bs\bs
\vbox{
\centerline{\bf 2. Proof of the main theorem}
\bsn
In this section we prove Theorem 1 using the estimate of the Castelnuovo-Mumford regularity. Since the notation of \cite{nwh} is used without further explanations, the reader is advised to read first few pages of the introduction as well as Sections 3.1 and 3.5 in \cite{nwh} before reading this section.}
\msn
{\bf 2.1.~Consequences of Corollary\,\,(1.4).} We first deduce upper estimates of the supports of $\Diff^2(\mu'''_{\rm max})$, $\Diff(\mu'')$, $\mu'$, $\mu'''_{{\rm def}}$ and lower estimates for $\Diff^2(\rho)$, $\Diff(\nu'')$, $\nu'''_{{\rm def}}$, $\nu'$ from Corollary~(1.4) using the self-duality \cite[Corollary~3.6]{nwh}.
\msn
{\bf Proposition~2.1.} {\it In the notation of \cite[Section 3.5]{nwh}, we have}
$$\Diff^2(\mu'''_{\rm max})_k=\Diff(\mu'')_k=\mu'_k=0\,\,\,(k>2d{-}2),\q\mu'''_{{\rm def},k}=0\,\,\,(k>2d{-}3).
\leqno(2.1.1)$$
\msn
{\it Proof.} We have the spectral sequence associated with the filtration $G$ in \cite[Section 3.1]{nwh}
$$E_1^{p,q}=H^{p+q}_{\mm}\Gr_G^pM\Longrightarrow H^{p+q}_{\mm}M,
\leqno(2.1.2)$$
where $\mm:=(x_1,\dots,x_4)\subset R$, and $\Gr_G^0M=M'''$, $\Gr_G^1M=M''$, $\Gr_G^2M=M'$. We can use generic projections $\pi_i:\C^4\to\C^i$ to calculate the local cohomology, see also the proof of \cite[Corollary 4.8]{Ei}. We then get (using the short exact sequence (2.3.8) below for $M'''$)
$$E_1^{p,q}=0\q\h{unless}\q(p,q)\in\bl\{(0,1),(0,2),(1,0),(2,-2)\br\}.
\leqno(2.1.3)$$
So the spectral sequence degenerates at $E_1$, and the regularity of $M'$, $M''$, $M'''$ is at most $2d{-}2$ by Corollary~(1.4) and Remark~(1.1)(ii). This implies (2.1.1) except for the assertion about $\mu'''_{{\rm def},k}$. As for the latter, the short exact sequence (2.3.8) below implies the isomorphisms
$$H^1_{\mm}M'''=H^0_{\mm}M'''_{\rm def}\,,\q H^2_{\mm}M'''=H^2_{\mm}M'''_{\rm max}\,,
\leqno(2.1.4)$$
where $H^j_{\mm}M'''=0$ ($j\ne 1,2$), $H^j_{\mm}M''''_{\rm max}=0$ ($j\ne 2$), $H^j_{\mm}M''''_{\rm def}=0$ ($j\ne 0$). Proposition~(2.1) then follows from Remark~(1.1)(ii).
\ms
By \cite[Corollary~3.6]{nwh}, we can deduce the following from Proposition~(2.1).
\msn
{\bf Corollary~2.1.} {\it In the notation of \cite[Section 3.5]{nwh}, we have}
$$\aligned\Diff^2(\rho)_k=0\,\,\,(k<2d{+}4),\q&\Diff(\nu'')_k=0\,\,\,(k<2d{+}3),\\ \nu'''_{{\rm def},k}=0\,\,\,(k<2d{+}2),\q&\nu'_k=0\,\,\,(k<2d{+}3).\endaligned
\leqno(2.1.5)$$
\msn
{\bf 2.2. Other estimates.} For the proof of Theorem~1, we also need lower estimates of the supports of $\mu''$, $\mu'$ and upper estimates for $\Diff(\nu'')$, $\nu'''_{\rm def}$, $\Diff(\nu'''_{\rm max})$ as follows.
\msn
{\bf Proposition~2.2.} {\it In the notation of \cite[Section 3.5]{nwh}, we have}
$$\Diff(\mu'')_k=\mu''_k=\mu'_k=0\q(k<d+3).
\leqno(2.2.1)$$
\msn
{\it Proof.} Let $J\subset R$ be the saturation of the Jacobian ideal $(\dd f)\subset R$ for the product of the associated primes of $(\dd f)$ with codimension 3 in $\C^4$ so that $R/J=M'''(4)$ in the notation of Section 3.1, {\lc} Then Proposition~(2.2) is reduced to Lemma~(2.2) below, since the latter implies that $M'_k=M''_k=0$ ($k<d{+}3$). (Here it is quite difficult to get an assertion about $\mu'''_{\rm def}$.).
\msn
{\bf Lemma~2.2.} {\it In the above notation, we have}
$$J\subset\mm^{d-1},\q\h{that is,}\q J_k=0\,\,\,(k<d{-}1).
\leqno(2.2.2)$$
\msn
{\it Proof.} We can define the ideal $J\subset R$ also by the condition that $g^{\sim}|_U\in\J_f(k)|_U$ for $g\in R_k$. Here $g^{\sim}$ is a section of $\OO_{\PP^3}(k)$ defined by $g$, $\J_f\subset\OO_{\PP^3}$ is the ideal generated by the Jacobian ideal $(\dd f)\subset R$, and $U\subset\PP^3$ is the complement of some finite subset of $\PP^3$. Taking a general hyperplane section, the assertion (2.2.2) is then reduced to the case $n=3$, where it is more or less known. Indeed, the assertion for $n=3$ is equivalent to that $\mu'_k=0$ for $k<d+2$ (since $\deg\dd f/\dd x_i=d{-}1$), or equivalently, for $k>2d{-}2$ by the symmetry of the $\mu'_k$, see \cite{DiSa}. But the latter follows from (1.4.4) for $n=3$, where we have
$$\nu_k=0\,\,\,(k<d+2),\q\mu''_k+\nu_{3d-k}=\tau\,\,\,(k\in\Z),$$
see \h{\lc} and also \cite{DIM}. So Lemma~(2.2) and Proposition~(2.2) follow.
\msn
{\bf Corollary~2.2.} {\it In the notation of \cite[Section 3.5]{nwh}, we have}
$$\aligned&\Diff(\nu'')_k=0\,\,\,(k>3d{-}2),\q\nu'''_{{\rm def},k}=0\,\,\,(k>3d{-}3),\\&\Diff(\nu'''_{\rm max})_k=0\,\,\,(k\notin[d{+}3,3d{-}1]).\endaligned
\leqno(2.2.3)$$
\msn
{\it Proof.} By \cite[Corollary~3.6]{nwh}, the first two inequalities follows from Proposition~(2.2), and the last one is reduced to the vanishing of $\Diff(\nu'''_{\rm max})_k$ for $k<d{+}3$, which follows from that of $\Om^3_k$ ($k<3$), since $\nu'''_{{\rm def},k}=0$ ($k<2d{+}2$) by (2.1.1) using {\lc} This finishes the proof of Corollary~(2.2).
\msn
{\bf 2.3.~Proof of Theorem~1.} Set
$$\aligned&\muc_k:=\mu_k\,,\q\nuc_k:=\nu_k-\nu'_k\,,\q\rhoc_k:=\rho_k-\nu'_k\,,\\&\chic_{f,k}:=\muc_k-\nuc_{k+d}+\rhoc_{k+2d}\q(k\in\Z).\endaligned$$
Define $\Euc_i^{\les m}$ in the same way as ${\rm Eu}_i^{\les m}$ in \cite[(4.2.9)]{nwh} by replacing $\chi_f$ with $\chic_f$. Put
$\Euc_i:=\Euc_i^{\les m}$ ($m\gg 0$) as in {\lc} We have by definition
$$\Euc_i^{\les m}={\rm Eu}_i^{\les m}\,\,\,(m\gg 0),\q\h{that is,}\q\Euc_i={\rm Eu}_i.$$
\sk
From Proposition~(2.1) and Corollary~(2.2) together with the relation $\mu_k-\nu_k+\rho_k=0$ ($k>4d-4$), we can deduce that 
$$\Diff(\muc)_k=\tfrac{1}{2}\1\Diff(\nuc)_{k+d}=\Diff(\rhoc)_{k+2d}=\tau_{Z'}\q(k\ges 2d),
\leqno(2.3.1)$$
where $\tau_{Z'}$ is the global Tjurina number of a general hyperplane section $Z'$ of $Z$. We thus get
$$\Diff(\chic_f)_k=0\q(k\ges 2d).
\leqno(2.3.2)$$
Since $\chic_{f,k}=\chi_{f,k}=0$ ($k\gg 0$), this implies the vanishing
$$\chic_{f,k}=0\q(k\ges 2d-1),
\leqno(2.3.3)$$
and then the inequality
$${\rm Eu}_d+1=\Euc_d+1=\chic_{f,d}=\chi_{f,d}-\nu'_{3d}\les\chi_{f,d}=-\chi(U)+1.
\leqno(2.3.4)$$
using Corollary~(2.1) (which implies that $\nu'_{2d}=0$) and also \cite[Corollary 4.8]{nwh} (which proves the last equality of (2.3.4), where \cite[Corollary 6.3]{WiYu} is used in an essential way).
\sk
For $m\gg 0$ and $i=d$, we thus get the following inequalities similar to \cite[(4.3.6)]{nwh}
$$-\chi(U)\les{}^{(\infty)}{\rm Eu}_i^{\les m}\les{\rm Eu}_i^{\les m}={\rm Eu}_i\les-\chi(U),
\leqno(2.3.5)$$
where ${}^{(\infty)}{\rm Eu}_i^{\les m}$ is as in {\lc} (The condition $m\gg 0$ is used for the equality ${\rm Eu}_d^{\les m}={\rm Eu}_d$.) All the inequalities in (2.3.5) then become equalities, and we get in particular
$$\nu'_{3d}=0,\q\h{that is,}\q\mu'''_{{\rm def},d}=0,
\leqno(2.3.6)$$
using \cite[Corollary 3.6]{nwh}. By Lemma~(2.3) below together with (2.3.1), this implies that
$$\aligned&\nu'_k=0\,\,\,(k\ges 3d),\q\mu'''_{{\rm def},k}=0\,\,\,(k\les d),\\&\tfrac{1}{2}\1\Diff(\nu)_{k+d}=\Diff(\rho)_{k+2d}=\tau_{Z'}\,\,\,(k>2d),\q\chi_{f,k}=0\,\,\,(k\ges 2d).\endaligned
\leqno(2.3.7)$$
Theorem~1 now follows from \cite[Theorems 4.3]{nwh} for $m=2d-1$ using Theorem~4.9, \lc, since $\mu_k-\rho_k$ is independent of $k\gg 0$. We have $\rho_{4d-1}^{(2)}=0$ by \cite[Proposition 1]{nwh}, and $\mu_{2d-1}^{(2)}=0$, since $\nu_k^{(2)}=0$ ($k\ges 3d$) and $\mu_{2d-1}^{(\infty)}=0$ by (15), {\lc} This finishes the proof of Theorem 1.
\msn
{\bf Lemma~2.3.} {\it In the above notation, let $k$ be an integer at most $d$. Then $\mu'''_{{\rm def},k}=0$ if $\mu'''_{{\rm def},k+1}=0$.}
\msn
{\it Proof.} We have the short exact sequence of graded $R$-modules as in \cite[(3.1.8)]{nwh}
$$0\to M'''\buildrel{\iota}\over\to M'''_{\rm max}\buildrel{\phi}\over\to M'''_{\rm def}\to 0.
\leqno(2.3.8)$$
Take any $v\in M'''_{{\rm def},k}$. Let $v'\in M'''_{{\rm max},k}$ such that $\phi(v')=v$. Let $x,y$ be the pull-backs of the coordinates of $\C^2$ by the generic projection $\pi_2:\C^4\to\C^2$ in \cite[Section 3.1]{nwh}. Since $xv=yv=0$, there are $v_1,v_2\in M'''_{k+1}$ such that
$$\iota(v_1)=xv',\q\iota(v_2)=yv'.$$
We have
$$yv_1=xv_2\q\h{in}\q M'''_{k+2},$$
since $\iota(yv_1)=\iota(xv_2)$. By Lemma~(2.2) there are canonical isomorphisms
$$M'''_j=R_{j-4}\q(j\les k+2).$$
(Here the assumption $k\les d$ is used.) We then get $v_0\in M'''_k$ such that
$$xv_0=v_1,\q yv_0=v_2,$$
since $R$ is a unique factorization domain. Moreover $\iota(v_0)=v'$ (hence $v=0$), since $M'''_{\rm max}$ is a free $\C[x,y]$-module (in particular, it is $x$-torsion-free). So Lemma~(2.3) follows.

\end{document}